\documentclass[onefignum,onetabnum,onealgnum,oneeqnum]{siamart190516}
\usepackage[papersize={8.5in,11in},margin=1.0in,footskip=0.3in]{geometry}
\usepackage{amssymb}
\usepackage{amsmath}
\usepackage{enumitem}
\usepackage{microtype}

\setlist[enumerate]{leftmargin=2.6em,label=(\roman*),topsep=0.5em,parsep=0.25em}
\setlist[itemize]{leftmargin=2em,topsep=0.25em,,parsep=0.25em}

\definecolor{bluey}{HTML}{1f80ad}
\hypersetup{urlcolor=bluey,linkcolor=bluey,citecolor=bluey}

\usepackage[font={sf,small},labelsep=period,labelfont=bf,margin=0.25in]{caption}

\definecolor{bluesection}{HTML}{175280}
\makeatletter
\renewcommand\section{\@startsection{section}{1}{.25in}%
                                   {1.3ex \@plus .5ex \@minus .2ex}%
                                   {-.5em \@plus -.1em}%
                                   {\reset@font\normalsize\bfseries\color{bluesection}}}
\renewcommand\subsection{\@startsection{subsection}{2}{.25in}%
                                     {1.3ex\@plus .5ex \@minus .2ex}%
                                     {-.5em \@plus -.1em}%
                                     {\reset@font\normalsize\bfseries\color{bluesection}}}
\makeatother

\crefname{figure}{Fig.\!}{Figs.\!}
\crefname{section}{\S\!}{\S\!}
\crefname{subsection}{\S\!}{\S\!}
\newcommand{\R}{\mathbb R}
\DeclareMathOperator{\sign}{sign}
\DeclareMathOperator{\range}{range}
\newcommand{\shorteq}{\resizebox{0.75\width}{\height}{\,=}}

\makeatletter
\newcommand\incircbin
{%
  \mathpalette\@incircbin
}
\newcommand\@incircbin[2]
{%
  \mathbin%
  {%
    \ooalign{\hidewidth$#1#2$\hidewidth\crcr$#1\bigcirc$}%
  }%
}
\newcommand{\cplus}{\incircbin{+}}
\newcommand{\cminus}{\incircbin{-}}
\newcommand{\circdot}{\incircbin{\cdot}}
\makeatother
						
\usepackage{algorithm}
\usepackage[noend]{algpseudocode}
\makeatletter
\algrenewcommand\ALG@beginalgorithmic{\footnotesize}
\renewcommand{\ALG@name}{\small Algorithm}
\makeatother

\usepackage{etoolbox}
\makeatletter
\newif\if@algorithm
\patchcmd{\@float@Hx}%
  {\@nodocument}%
  {\@nodocument%
   \ifnum\pdfstrcmp{#1}{algorithm}=0 %
     \@algorithmtrue%
     \setlength{\columnwidth}{\algorithmwidth}%
   \fi}%
  {}{}%
\renewcommand\float@makebox[1]{%
  \hbox{%
    \if@algorithm\hspace*{-.5\dimexpr\algorithmwidth-\textwidth}\fi%
    \vbox{\hsize=#1 \@parboxrestore
      \@fs@pre\@fs@iftopcapt
        \ifvoid\@floatcapt\else\unvbox\@floatcapt\par\@fs@mid\fi
        \unvbox\@currbox
      \else\unvbox\@currbox
        \ifvoid\@floatcapt\else\par\@fs@mid\unvbox\@floatcapt\fi
      \fi\par\@fs@post\vskip\z@}}}
\makeatother
\newlength{\algorithmwidth}
\AtBeginDocument{\setlength{\algorithmwidth}{0.9\textwidth}}

\title{A Connected Component Labeling Algorithm for Implicitly-Defined Domains}
\author{Robert I.~Saye\thanks{Lawrence Berkeley National Laboratory, Berkeley, California, USA (\texttt{rsaye@lbl.gov})}}
\date{\today}

\begin{document}

\maketitle

\begin{abstract}
A connected component labeling algorithm is developed for implicitly-defined domains specified by multivariate polynomials. The algorithm operates by recursively subdividing the constraint domain into hyperrectangular subcells until the topology thereon is sufficiently simple; in particular, we devise a topology test using properties of Bernstein polynomials. In many cases the algorithm produces a certificate guaranteeing its correctness, i.e., two points yield the same label if and only if they are path-connected. To robustly handle various kinds of edge cases, the algorithm may assign identical labels to distinct components, but only when they are exactly or nearly touching, relative to a user-controlled length scale. A variety of numerical experiments assess the effectiveness of the overall approach, including statistical analyses on randomly generated multi-component geometry in 2D and 3D, as well as specific examples involving cusps, self-intersections, junctions, and other kinds of singularities.
\end{abstract}



\section{Introduction}

In this paper, we develop a connected component labeling algorithm for implicitly-defined domains specified by multivariate polynomials. In particular, we consider domains of the form $\Omega := \mathcal U \setminus \{\phi = 0\}$, where $\phi: \R^d \to \R$ is a given polynomial and $\mathcal U \subset \R^d$ is a given bounded $d$-dimensional hyperrectangle, and consider the problem of implementing a labeling function $\chi : \Omega \to \mathbb N$ such that $\chi(x) = \chi(y)$ if and only if the two points $x, y \in \Omega$ are path-connected.

Our interest in this labeling problem stems from the author's prior work on quadrature algorithms for multi-component implicitly-defined domains \cite{SayeSP}. For a general, piecewise-smooth integrand function $f$, these algorithms output a quadrature scheme of the form $\int_{\mathcal U} f \approx \sum_i w_i f(x_i)$ and have the following key property: if $V$ is a connected component of $\mathcal U \setminus \{\phi = 0\}$ and $f$ is sufficiently smooth on $V$, then $\int_V f \approx \sum_{x_i \in V} w_i f(x_i)$ represents a high-order accurate quadrature scheme on $V$. In other words, to build a quadrature scheme on $V$, one may simply discard the quadrature points outside of $V$, leaving the weights of the remaining points unmodified. Consequently, one may think of the bulk quadrature scheme over $\mathcal U$ as an agglomeration of smaller quadrature schemes over the individual connected components of $\mathcal U \setminus \{\phi = 0\}$. The algorithms developed in \cite{SayeSP}, however, do not automatically determine this grouping; one of the aims of the present work is to design a simple and efficient algorithm for this purpose. Besides this application, connected component analysis arises in a wide variety of different settings, including in computer-aided geometric design, shape modeling and pattern matching, and in building topological feature descriptors of physical data sets \cite{10.1007/978-3-319-04099-8_6}.

A number of different approaches could be taken to solve the connected component labeling problem:
\begin{itemize}
\item One approach is to apply the workhorse tools of real algebraic geometry, such as the \textit{cylindrical algebraic decomposition} (CAD) algorithm pioneered by Collins \cite{10.1007/3-540-07407-4_17,caviness2012quantifier}, or by computing roadmaps \cite{10.1145/281508.281533,10.1007/3-540-33099-2}. In exact or arbitrary-precision arithmetic, CAD and roadmap methods can be used to identify and analyze the connected components of general, arbitrarily complex semi-algebraic sets (of which $\Omega$ is one particular instance). However, these methods can be computationally expensive; moreover, their implementation in fixed-precision arithmetic requires considerable care, even for relatively simple geometry.
\item A second approach is to apply the computational methods of Morse theory, Reeb graphs, and contour trees \cite{matsumoto2002introduction,10.1007/978-3-319-04099-8_6}. Roughly speaking, these methods compute the topology of the set of level sets of $\phi$ via the location of its extrema and saddle points, combined with gradient path tracing algorithms. These methods become increasingly more intricate to implement as the dimension $d$ increases. 
\item A third approach might be to isolate individual components of $\Omega^+ := \mathcal U \cap \{\phi > 0\}$ by finding a polynomial which can separate them (followed by a similar procedure for $\Omega^- := \mathcal U \cap \{\phi < 0\}$). For example, if we could find a polynomial $u: \R^d \to \R$ such that  $\Omega^+ \cap \{ u = 0 \}$ is empty and $u(x)u(y) < 0$ for two sample points $x,y \in \Omega^+$, it would follow that the zero level set of $u$ forms a kind of divider that separates the two connected components associated with $x$ and $y$. To find this polynomial, one could use polynomial positivstellensatz or sums-of-squares methods \cite{POWERS199899,putinar2008emerging} to compute, if possible, a polynomial $w$ such that: (i) $w > 0$ on $\Omega^+$; (ii) $w = u^2$; and (iii) $u(x)u(y) < 0$. Part (i), (ignoring conditions (ii) and (iii) and various other subtleties,) can be solved efficiently via linear or semidefinite programming techniques, see, e.g., \cite{putinar2008emerging,doi:10.1137/18M118935X}; however, conditions (ii) and (iii) end up creating a quadratically-constrained quadratic program for the coefficients of $w$, whose associated feasibility condition involves a challenging non-convex constraint. These ideas, among various others for computing separating polynomials, were investigated as part of the present work; however, no sufficiently elegant or efficient approach was found.
\item Finally, a fourth and considerably simpler approach is to subdivide the constraint domain $\mathcal U$ into small enough cells for which the topology of $\Omega$ is easy to determine. Depending on the tessellation method, this approach can be quite effective (fast, simple to implement, high resolution power), but in some cases the subdivision process may need to be stopped, potentially leading to some components being incorrectly merged or broken.
\end{itemize}
The preceding discussion serves to highlight that, depending on the final application, one must ultimately decide on a suitable compromise between: (i) absolute correctness/certifiability of the connected component labeling algorithm, potentially needing costly arbitrary-precision computation (as exemplified by CAD-based methods); (ii) robustness, e.g., in suitably handling fixed-precision roundoff errors or uncertainty in the input polynomial coefficients; (iii) computational complexity of constructing the labeling function as well as its subsequent evaluation; and (iv) implementation simplicity. The target application of this work concerns the development of numerical methods for multi-physics simulations involving highly complex geometry such as liquid atomization dynamics; in particular, the input polynomial $\phi$ represents the fluid interface and is computed dynamically. This application necessitates prioritizing aspects (ii), (iii), and (iv), which, in turn, necessitates a design choice for how to handle edge cases such as nearly touching components or interfacial self-intersections. Our choice here is to require that components are never broken, i.e., if $x,y$ are path-connected in $\Omega$ and $\chi(x), \chi(y)$ is the output of the labeling algorithm, then $\chi(x) = \chi(y)$ is an absolute certainty. It follows that the labeling algorithm must be permitted to merge or ``glue'' together distinct components in the (presumably rare) cases of uncertain topology. The algorithm developed in this work glues components only if they are exactly or nearly touching, relative to a user-defined gap threshold, typically orders of magnitude smaller than the length scale of $\mathcal U$.

With this design objective in mind, the connected component labeling algorithm developed here follows the fourth approach mentioned earlier: it operates by recursively subdividing the constraint hyperrectangle $\mathcal U \subset \R^d$ into hyperrectangular subcells until the topology of $\Omega$ thereon is sufficiently simple or has reached a smallest-permitted size. This is achieved through a combination of quadtrees (in 2D) or octrees (in 3D), along with simple yet effective topology tests using properties of Bernstein polynomials. Once the tree is constructed, its leaf cells form the vertices of a graph whose edges are determined by sign attainability tests on the faces of adjacent leaf cells. The connected components of this graph ultimately decide the overall labeling: $\chi(x)$ equals the label of the leaf cell containing $x$, the latter found by fast tree traversal methods. We show in this paper that the overall approach is efficient in handling various degrees of geometric complexity. In particular, if the recursive subdivision process terminates without reaching the user-defined smallest subcell size, then the labeling is provably correct. Applied to a fixed class of randomly generated geometry, the percentage of cases in which the algorithm is uncertain about the topology decreases exponentially as a function of the maximum tree depth; moreover, only a small fraction of these cases yield artificially-glued components.

An outline for the rest of the paper is as follows. In \cref{sec:prelim}, we establish some preliminaries on Bernstein polynomials and define the simply connected topology test. \cref{sec:algo} presents the main connected component labeling algorithm, followed by a discussion of its features. Numerical tests are presented in \cref{sec:results}, analyzing the algorithm's success and failure rates on randomly generated multi-component geometry in 2D and 3D as well as its behavior on particular examples exhibiting cusps, self-intersections, and other kinds of singularities. Concluding remarks with a brief discussion of possible extensions are given in \cref{sec:conclusion}.

\section{Preliminaries}
\label{sec:prelim}

A key part of the connected component labeling algorithm is an effective means for evaluating the range of attainable values of a polynomial. One of the most accurate and straightforward methods for doing so is through the use of the Bernstein basis \cite{garloff1993bernstein,Rajyaguru2017,lin1995methods,beska1989convexity}. These methods make use of a convex hull property and guarantee that a polynomial's value, at any point in its rectangular reference domain, is no larger (or smaller) than its maximum (or minimum) coefficient in the Bernstein basis. Especially useful in the present setting, these bounds become monotonically more accurate under the stable operations of Bernstein subdivision. In addition, owing to the hyperrectangular constraint domain $\mathcal U$ and its subdivision into hyperrectangular subcells, it is particularly natural to use a tensor-product basis. Accordingly, a tensor-product Bernstein basis is adopted throughout this work.

In $d$ dimensions, let $\phi$ be a tensor-product Bernstein polynomial of degree $n = (n_1,\ldots,n_d)$, defined relative to the hyperrectangle $U = [\alpha_1, \beta_1] \times \cdots \times [\alpha_d, \beta_d]$; $\phi$ takes the form
\[ \phi(x_1, \ldots, x_d) = \sum_{i_1 = 0}^{n_1} \cdots \sum_{i_d = 0}^{n_d} c_{i_1,\ldots,i_d} b_{i_1}^{n_1,[\alpha_1,\beta_1]}(x_1) \cdots b_{i_d}^{n_d,[\alpha_d,\beta_d]}(x_d) = \sum_{i \in {\mathbb N}^n} c_i b_i^{n,U}(x) \]
where
\[ b_\ell^{\eta,[\alpha,\beta]}(x) = \frac{\binom{\eta}{\ell}}{(\beta - \alpha)^\eta} (\beta - x)^{\eta-\ell} (x - \alpha)^\ell, \quad \ell = 0, \ldots, \eta,\]
are the one-dimensional Bernstein basis functions of degree $\eta$ relative to the interval $[\alpha,\beta]$, $c_i = c_{i_1,\ldots,i_d}$ denotes the $i$th Bernstein coefficient of $\phi$ for a multi-index $i \in {\mathbb N}^n := [0,n_1] \times \cdots \times [0,n_d]$, and $\smash{b_i^{n,U}}(x) = \smash{\prod_{j=1}^d} \smash{b_{i_j}^{n_j,[\alpha_j,\beta_j]} (x_{i_j})}$ denotes the product of basis functions. This representation gives the Bernstein coefficients of $\phi$ relative to the given hyperrectangle $U$; subdivision refers to transforming these coefficients relative to a subset hyperrectangle and can be stably computed using the de Casteljau algorithm \cite{FAROUKI19881}.

In the connected component labeling algorithm, we need a method to assess whether the topology of $U \cap \{\phi > 0\}$ and $U \cap \{\phi < 0\}$ is simple enough. This is achieved through the following three concepts.\footnote{\cref{def:cmi} has, in various guises, appeared before in the literature; the remaining two concepts (\cref{def:sea,def:sc}, along with their implications) are perhaps new, and target specifically the objectives of the present work.}

\begin{definition} \label{def:cmi}
A Bernstein polynomial $\phi = \sum_{i \in {\mathbb N}^n} c_i b_i^{n,U}(x)$ is called \emph{coefficient monotone increasing on $U$ (in the direction $k$)} if there exists a coordinate direction $k$ such that $c_i \leq c_{i + e_k}$ for all $i \in {\mathbb N}^n$, $i_k < n_k$; here, $e_k$ denotes the standard basis vector in the direction of the $k$th coordinate. The polynomial is called \emph{coefficient monotone} if $\phi$ or $-\phi$ is coefficient monotone increasing.
\end{definition}

An equivalent viewpoint comes from a convenient property of differentiation in the Bernstein basis: up to a multiplicative factor, the coefficients of the derivative are formed via first-order divided differences of the input polynomial's coefficients. Therefore, $\phi$ is coefficient monotone if and only if the Bernstein coefficients of its derivative in the corresponding direction are either all non-negative or all non-positive. It follows that coefficient monotone polynomials are monotone in the conventional sense and therefore attain their extrema on the corresponding faces of the hyperrectangle. Tests of Bernstein coefficient monotonicity have a variety of applications, including, e.g., in polynomial range evaluation \cite{10.1007/s10898-008-9382-y}. In the present application, an important property is that, for a coefficient monotone polynomial $\phi$, for any $x \in U \cap \{\phi > 0\}$, $x$ can be path-connected to one of the corresponding faces of the hyperrectangle without leaving the region $U \cap \{\phi > 0\}$, and analogously for points in $U \cap \{\phi < 0\}$.

We next establish a sufficiently accurate means to determine whether a polynomial attains positive or negative values on a given hyperrectangle. The following definition sets the requirements on an algorithm implementing this task.

\begin{definition} \label{def:sea}
Given a Bernstein polynomial $\phi = \sum_{i \in {\mathbb N}^n} c_i b_i^{n,U}(x)$ on the hyperrectangle $U$, a $d$-dimensional \emph{sign evaluation algorithm} $\sigma_d(\phi, U)$ outputs the possible signs of $\phi$ on $U$ and must satisfy the following properties:
\begin{enumerate}
\item $\sigma_d(\phi,U) \subseteq \{-1, 0, +1\}$;
\item $\sigma_d(\phi,U)$ must be a superset of the ground-truth, i.e., $\{\sign(\phi(x)) : x \in U\} \subseteq \sigma_d(\phi,U)$, where $\sign(u)$ is the standard sign operator,
\[ \sign(u) = \left\{\begin{array}{rl} +1, & \text{if $u > 0$,} \\ 0, & \text{if $u = 0$,} \\ -1, & \text{if $u < 0$.} \end{array} \right. \]
\item if $c_i \geq 0$ for all $i$, then $-1 \notin \sigma_d(\phi,U)$; if $c_i \leq 0$ for all $i$, then $+1 \notin \sigma_d(\phi,U)$;
\item if $\phi$ is coefficient monotone in the direction $k$, then $\sigma_d(\phi,U) \subseteq \sigma_{d-1}(\phi_\alpha,U_\alpha) \cup_0 \sigma_{d-1}(\phi_\beta,U_\beta)$ must hold, where $\phi_\alpha$ and $\phi_\beta$ denote the restriction of $\phi$ to the corresponding lower $U_\alpha$ and upper $U_\beta$ faces of $U$, and $\sigma_{d-1}$ is a $(d-1)$-dimensional sign evaluation algorithm;\footnote{Here, $\cup_0$ denotes the union operation which includes zero when appropriate: if $s_1 \cup s_2 = \{-1,+1\}$, then $s_1 \cup_0 s_2 := \{-1,0,+1\}$; if $s_1 \cup s_2 \neq \{-1,+1\}$, then $s_1 \cup_0 s_2 := s_1 \cup s_2$.}
\item finally, a zero-dimensional sign evaluation algorithm must satisfy $\sigma_0(c,U) = \sign(c)$.
\end{enumerate}
\end{definition}

Note that a method which exactly computes the minimum and maximum value of $\phi$ on $U$ would trivially yield a sign evaluation algorithm. However, computing the extrema of arbitrary polynomials is a non-trivial and computationally expensive task, perhaps as difficult as the connected component labeling problem itself. The conditions of a sign evaluation algorithm are weaker than an exact computation, thereby allowing for simpler and more efficient algorithms; nevertheless, some degree of certitude is required. In particular, condition (iii) requires an implementation be at least as accurate as conventional Bernstein range evaluation; this ensures it inherits at least the same level of accuracy under the action of subdivision. Condition (iv) says that if $\phi$ happens to be coefficient monotone, the sign evaluation must be at least as accurate as what can be determined solely from the corresponding faces where $\phi$'s extrema are known to occur. One particularly important property is that if a sign evaluation algorithm outputs either $\{-1\}$ or $\{+1\}$, then we may conclude, with absolute certainty, $\phi$ is uniformly signed throughout $U$. An example of a simple sign evaluation algorithm is given in the next section.

Finally, a notion of whether the topology of $U \setminus \{\phi = 0\}$ is ``simple enough'' is recursively defined as follows. This definition implicitly depends on the presence of a suitable sign evaluation algorithm $\sigma_d$, so it is assumed that one has been implemented and fixed ahead of time.

\begin{definition} \label{def:sc}
A $(d \geq 1)$-dimensional Bernstein polynomial $\phi = \sum_{i \in {\mathbb N}^n} c_i b_i^{n,U}(x)$ is called \emph{simply connected on $U$} if at least one of the following three conditions hold: (i) $\sigma_d(\phi,U) = \{+1\}$; or (ii) $\sigma_d(\phi,U) = \{-1\}$; or (iii) there is coordinate direction $k$ such that (a) $\phi$ is coefficient monotone in the direction $k$ \textsc{and} (b) $\phi$ is simply connected on the corresponding upper and lower faces of $U$. By definition, a zero-dimensional Bernstein polynomial is \emph{simply connected}.
\end{definition}

A simply connected polynomial guarantees a simple topology of both $\Omega^- := U \cap \{\phi < 0\}$ and $\Omega^+ := U \cap \{\phi > 0\}$, as follows. If conditions (i) or (ii) in \cref{def:sc} are met, then $\phi$ is non-zero throughout the hyperrectangle and so exactly one of $\Omega^\pm$ is empty and the other the whole hyperrectangle. Otherwise, there is a coordinate axis $k$ such that any point in $x \in \Omega^\pm$ can be path-connected to one of the corresponding faces. These faces are themselves simply connected; inductively it follows that if $\Omega^+$ (resp., $\Omega^-$) is nonempty, then $\Omega^+$ (resp., $\Omega^-$) has exactly one connected component.\footnote{In fact, the sets $\Omega^+$ or $\Omega^-$ are simply connected in the conventional topological sense, but it is unnecessary to prove this fact in the present work.}

To prove the correctness of the connected component labeling algorithm, it is useful to establish the following two properties of simply connected polynomials:
\begin{itemize}
\item The first property is that if $\phi$ is simply connected on a hyperrectangle $U$, then $\phi$ is simply connected on every face of $U$. This can be shown via induction on the dimension; we illustrate here with a three-dimensional example. Suppose $\phi$ is simply connected on a rectangular 3D prism. If $\sigma_d(\phi) \in \bigl\{ \{-1\}, \{+1\} \bigr\}$, then it is uniformly signed and thus trivially simply connected on every face of the prism. Otherwise, $\phi$ is coefficient monotone in the up direction, say, such that the restriction of $\phi$ to the bottom and top faces is simply connected; in particular, by induction, we have that $\phi$ is simply connected on every edge of the bottom and top face. Now, on each vertical face of the prism (those faces excluding the bottom and top), the restriction of $\phi$ is coefficient monotone in the up direction, and, as was just concluded, $\phi$ is simply connected on the corresponding lower and upper edges. Therefore, $\phi$ is simply connected on every vertical face. The preceding arguments can be extended to any dimension, and the one-dimensional base case of the inductive argument trivially holds.

\item The second property is that if $\phi$ is simply connected on a hyperrectangle $U$, then a sign evaluation algorithm yields \textit{exact} results, i.e., the output of $\sigma_d(\phi,U)$ is precisely $\{\sign(\phi(x)): x \in U\}$. Similar to before, this can be shown via induction. Suppose $\phi$ is simply connected on $U$. If $\sigma_d(\phi,U) \in \smash{\bigl\{ \{-1\}, \{+1\} \bigr\}}$, then $\phi$ is with absolute certainty uniformly signed and $\sigma_d$ is exact. Otherwise, if $\phi$ is coefficient monotone on $U$ in the up direction, say, then $\sigma_d(\phi,U)$ is at least as accurate as $\sigma_{d-1}$ applied to the restriction of $\phi$ on the lower and upper faces of $U$; by induction, the latter calculation is exact, and therefore $\sigma_d$ is, too.
\end{itemize}
Combining these two properties, we observe that if $\phi$ is simply connected on a hyperrectangle $U$, then a sign evaluation algorithm yields \textit{exact} results on every face of $U$.

\section{Connected Component Labeling Algorithm}
\label{sec:algo}

\begin{algorithm}[t]
\caption{\small Given a $d$-dimensional, degree $n = (n_1,\ldots,n_d)$ Bernstein polynomial $\phi = \sum_{i \in {\mathbb N}^n} c_i b_i^{n,U}(x)$ and a corresponding hyperrectangular domain $U = [\alpha_1, \beta_1] \times \cdots \times [\alpha_d, \beta_d]$, evaluate $\range(\phi,U)$.}
\label{algo:range}
\begin{algorithmic}[1]
	\If{$d = 1$ and $n_1 \leq 4$}
		\State Compute $\displaystyle \theta := \min_{\alpha_1 \leq x \leq \beta_1} \phi(x)$ and $\displaystyle \Theta := \max_{\alpha_1 \leq x \leq \beta_1} \phi(x)$. \label{line:lowdegree}
	\Else
		\State Set $\theta := \min_i c_i$ and $\Theta := \max_i c_i$.
		\If{$d > 1$}
		\For{$k = 1, 2, \ldots, d$}
			\If{$\phi$ is coefficient monotone in direction $k$ on $U$}
				\State Compute $R := \range(\phi|_{x_k = \alpha_k}, \text{lower face of $U$}) \cup \range(\phi|_{x_k = \beta_k}, \text{upper face of $U$})$. \label{line:facerange}
				\State Update $\theta \leftarrow \max \{ \theta, \inf R \}$. \label{line:u1}
				\State Update $\Theta \leftarrow \min \{ \Theta, \sup R \}$.  \label{line:u2}
			\EndIf
		\EndFor
		\EndIf
	\EndIf
	\State \textbf{return} $[\theta, \Theta]$.
\end{algorithmic}
\end{algorithm}

With these preliminaries established, we are now in a position to describe the connected component labeling algorithm. First, we describe a sufficiently-accurate Bernstein polynomial range evaluation algorithm, which approximately (and in various cases, exactly) evaluates the range of $\phi$ on a given hyperrectangle $U$. Our particular implementation is given in \cref{algo:range} and has the following properties:
\begin{itemize}
\item The output is an interval such that $[\inf_{x \in U} \phi(x), \sup_{x \in U} \phi(x)] \subseteq \range(\phi,U)$ always holds.
\item The algorithm is at least as accurate as conventional Bernstein range evaluation, i.e., $\range(\phi,U) \subseteq [\min_i c_i, \max_i c_i]$, the latter interval bounding the minimum and maximum value of the input's Bernstein coefficients.
\item For the simple cases of linear, quadratic, or cubic polynomials in $d=1$ dimensions, it is a particularly simple and efficient task to exactly evaluate the range of $\phi$; this is implemented on line \ref{line:lowdegree} and is an easy tweak benefiting overall accuracy.
\item If $\phi$ happens to be coefficient monotone on $U$, then the range evaluation algorithm is as least as accurate as what can be determined by evaluating the range on the corresponding faces of $U$ where the extrema of $\phi$ is known to occur (line \ref{line:facerange}). Further, if $\phi$ happens to be coefficient monotone in multiple directions, then the tightest possible range is returned (lines \ref{line:u1}--\ref{line:u2}).
\end{itemize}

The main purpose of the range evaluation algorithm is to implement a sign evaluation algorithm. In particular, we define
\[ \sigma_d(\phi,U) := \bigl\{ \sign(r) : r \in \range(\phi,U) \bigr\}. \]
In other words, $\sigma_d(\phi)$ is defined by the possible signs of $\phi$ as computed by the range evaluation algorithm. By construction, it satisfies all the required properties set out in \cref{def:sea}. In particular, it is exact for low-degree one-dimensional polynomials, whenever $0 \notin \range(\phi)$, as well as various other cases. Using the sign evaluation algorithm, it is also straightforward to implement an algorithm that tests whether $\phi$ is simply connected on a given hyperrectangle $U$; our implementation carries out the algorithmic process naturally suggested by \cref{def:sc}.

With these ingredients, we can now describe the connected component labeling algorithm. Given a Bernstein polynomial $\phi$ and a hyperrectangular constraint domain $\mathcal U$, we first test whether $\phi$ is simply connected on $\mathcal U$. If it is, then $\Omega^+ := \mathcal U \cap \{\phi > 0\}$ and $\Omega^- := \mathcal U \cap \{\phi < 0\}$ are either empty or have exactly one component, in which case the labeling problem is trivial. Otherwise, we recursively subdivide $\mathcal U$ into smaller and smaller subcells until either $\phi$ is simply connected on the subcell or a maximum recursion depth is met. In this work, we have opted for a quadtree/octree-style subdivision wherein the subcells maintain the aspect ratio of $\mathcal U$, though other subdivision algorithms are certainly possible. Upon conclusion of the subdivision process, we have a grid in which most (and in many cases, all) leaf cells have simple topology. In the second phase of the algorithm, two graphs ${\mathcal G}^\pm$ are constructed, each representing the connectivity of the cells overlapping $\Omega^+$ and $\Omega^-$. For example, if $\sigma_{d-1}$ applied to a shared face $\mathcal F = U_i \cap U_j$ of the grid yields $s \subseteq \{-1,0,+1\}$ and $-1 \in s$ (resp., $+1 \in s$), then an edge $(U_i, U_j)$ is inserted into ${\mathcal G}^-$ (resp., ${\mathcal G}^+$). These edges reflect the fact that there is likely (and in most cases, definitely) a path connecting $U_i \cap \Omega^\pm$ and $U_j \cap \Omega^\pm$. Ultimately, it is the connected components of ${\mathcal G}^-$ and ${\mathcal G}^+$ which define the possible labels of the overall algorithm; computing this labeling can be done with negligible cost via efficient disjoint-set/merge-find data structures. In essence, these graphs establish the globally topology of $\Omega$, whereas each leaf cell of the tree handles the local topology. Finally, to determine which component an arbitrary point $x \in \Omega$ belongs, we apply standard (and fast) tree traversal algorithms to find a leaf cell containing $x$, and then adopt its corresponding label from ${\mathcal G}^+$ (if $\phi(x) > 0$) or ${\mathcal G}^-$ (if $\phi(x) < 0$).

\begin{algorithm}[t]
\caption{\small Construction phase of the connected component labeling algorithm: given the polynomial $\phi$ and the hyperrectangular constraint domain $\mathcal U$, build the subdivision tree and create labels for its path-connected leaf cells.}
\label{algo:tree}
\begin{algorithmic}[1]
	\Statex \textit{\hspace{-2.5em} Construct the subdivision tree:}
	\State Create the root node $\nu_{\text{root}}$ representing $\mathcal U$.
	\State \textbf{call} \textsc{visit}($\nu_{\text{root}}$, $\mathcal U$, $0$).
	\Procedure{visit}{$\nu$, $U$, $\ell$}
		\If{$\phi$ is simply connected on $U$}
			\State Designate $\nu$ as a simply connected leaf node.
		\ElsIf{$\ell \geq \Lambda$} \label{line:maxdepth}
			\State Designate $\nu$ as a non-simply connected leaf node. \label{line:nonsc}
		\Else
			\State Designate $\nu$ as a non-leaf node, pointing to the following child nodes.
			\For{$i = 1, \ldots, 2^d$}
				\State Let $U_i$ denote the $i$th subcell of $U$.
				\State Insert a new node representing $U_i$ into the tree; call it $\nu_i$.
				\State \textbf{call} \textsc{visit}($\nu_i$, $U_i$, $\ell + 1$).
			\EndFor
		\EndIf
	\EndProcedure
	\Statex \textit{\hspace{-2.5em} Build subcell connectivity:}
	\State Initialize ${\mathcal G}^+$ and ${\mathcal G}^-$ as empty graphs.
	\For{every pair of leaf cells $U_i$, $U_j$ sharing a face}
		\State Let ${\cal F} := U_i \cap U_j$ denote the shared face. \label{line:face1}
		\State Compute $s := \sigma_{d-1}(\phi|_{\mathcal F}, \mathcal F)$.
		\If{$-1 \in s$}
			\State Insert the edge $(U_i, U_j)$ into ${\mathcal G}^-$. 
		\EndIf
		\If{$+1 \in s$}
			\State Insert the edge $(U_i, U_j)$ into ${\mathcal G}^+$. \label{line:face2}
		\EndIf
	\EndFor
	\Statex \textit{\hspace{-2.5em} Generate connected component labels:}
	\State Create a unique integer label for the combined set of connected components of ${\mathcal G}^-$ and ${\mathcal G}^+$.
	\State Initialize every leaf cell's $\cminus$ and $\cplus$ labels as \textit{undefined}.
	\For{every leaf cell $U$}
		\If{$U$ is a vertex of ${\mathcal G}^-$}
			\State Set the leaf cell's $\cminus$ label to the corresponding label of ${\mathcal G}^-$.
		\EndIf
		\If{$U$ is a vertex of ${\mathcal G}^+$}
			\State Set the leaf cell's $\cplus$ label to the corresponding label of ${\mathcal G}^+$.
		\EndIf
	\EndFor
\end{algorithmic}
\end{algorithm}

\begin{algorithm}[t]
\caption{\small Connected component label evaluation algorithm: using the subdivision tree constructed by \cref{algo:tree}, evaluate $\chi: \Omega \to \mathbb N$ for a given point $x \in \Omega := {\mathcal U} \setminus \{\phi = 0\}$.}
\label{algo:eval}
\begin{algorithmic}[1]
	\State Using standard, fast tree traversal methods, find a leaf cell $U$ containing $x$; if $x$ belongs to multiple leaf cells it does not matter which is used.
	\State Let $\circdot$ denote $\cminus$ if $\phi(x) < 0$ or $\cplus$ if $\phi(x) > 0$.
	\If{$U$'s $\circdot$ label is undefined}
		\State Create a new, unique label and assign it to $U$.
	\EndIf
	\State \textbf{return} the $\circdot$ label assigned to $U$.
\end{algorithmic}
\end{algorithm}

\cref{algo:tree} and \cref{algo:eval} sketch the implementation of the connected component labeling algorithm. A few complementary remarks are in order:
\begin{itemize}
\item In \cref{algo:tree}, whenever the restriction of $\phi$ to a subcell $U$ or its face is required, the corresponding Bernstein coefficients can be efficiently and accurately computed via the de Casteljau algorithm.
\item On line \ref{line:maxdepth}, $\Lambda$ specifies the user-defined maximum recursion depth. For example, if the constraint domain is the unit cube $[0,1]^d$, the parameter corresponds to the smallest-possible subcell having width $2^{-\Lambda}$. The effect of $\Lambda$ on the overall labeling is studied in detail in the next section.
\item \cref{algo:tree} only assigns labels to cells that are path-connected to other cells; as such, if a cell contains a part of $\Omega^+$, say, that has not been path-connected to any other cell overlapping with $\Omega^+$, then that cell does not (yet) have a $\cplus$ label defined. Our main motivation for this mainly concerns the scenarios for which the smallest-possible cell size is reached whereon the topology of $\Omega^+$ is still uncertain: in these cases, we do not want to unnecessarily create a new label when $\Omega^+$ thereon may actually be empty. In a sense, if a sign evaluation algorithm is unable to certify the presence of $\Omega^+$, then doing so will be up to the subsequent evaluation of $\chi$; indeed, \cref{algo:eval} creates a label only when a point $x \in \Omega^+$ is given which proves the nonemptiness of $\Omega^+$ in that cell.
\item By caching the results of the simply connected tests on the leaf cells, as computed by the first phase of \cref{algo:tree}, various steps in the second and last phase can be accelerated. For example, it is unnecessary to compute $s := \sigma_{d-1}(\phi|_{\mathcal F}, \mathcal F)$ on a shared face $\mathcal F$ if it is already known that $\phi$ is uniformly signed on one of the corresponding cells. On a related note, if a path connection already exists between $U_i$ and $U_j$, it is unnecessary to execute lines \ref{line:face1}--\ref{line:face2} because these lines would not alter that connectivity; thus, some faces can be skipped during the second phase of the algorithm.
\item Our particular implementation of the subdivision tree represents each node by a bit-field of width 16 bits. One bit specifies whether the node is a leaf or non-leaf node; the remaining 15 bits of a non-leaf node points to the children of that node; the remaining 15 bits of a leaf node store the type of the leaf cell (2 bits, flagging whether it is uniformly negative, uniformly positive, mixed sign and simply connected, or not simply connected) and 6 bits each to cache the $\cplus$ and $\cminus$ labels created by \cref{algo:tree,algo:eval}. This implies an upper bound of $2^6 - 1 = 63$ components for each phase, ample for our purposes. In addition, our implementation of \cref{algo:tree} bypasses the explicit creation of the graphs ${\mathcal G}^\pm$; instead, an efficient disjoint-set/merge-find data structure is used to represent the edge formations and subsequent label creation and lookup.
\item Finally, our implementation also makes use of a fuzzy threshold to robustly handle roundoff error in fixed-precision arithmetic. A tolerance $\epsilon$ is incorporated into the sign evaluation algorithm as follows: if $R =\range(\phi,U)$, then $\sigma_d(\phi,U) = \{+1\}$ iff $\inf R \geq \epsilon$, $\sigma_d(\phi,U) = \{-1\}$ iff $\sup R \leq -\epsilon$, and $\sigma_d(\phi) = \{-1,0,+1\}$ in all other cases. This tolerance helps to prevent inconsistencies arising from roundoff error, such as what might occur when the zero level set $\{\phi = 0\}$ grazes the boundary of a cell. Here, the tolerance is chosen to scale relative to the Bernstein coefficients of $\phi$ on $\mathcal U$; in particular, we have set $\epsilon = 10^3\, \epsilon_0 \max_i |c_i|$, where $\epsilon_0$ is machine epsilon (approximately $10^{-16}$ in double-precision arithmetic).
\end{itemize}

To conclude this section, we prove the correctness of the connected component labeling algorithm in the case the subdivision process finishes without hitting the maximum-imposed recursion depth. In other words, every leaf cell of the tree satisfies the simply connected property. In the prior section we saw that the sign evaluation algorithm is \textit{exact} in this case, and so two adjacent cells $U_i$ and $U_j$ are $\cplus$-connected if and only if $(U_i \cap U_j) \cap \{\phi > 0\}$ is nonempty. It follows that if $x, y \in \Omega^+ := {\mathcal U} \cap \{\phi > 0\}$ are path connected, then a corresponding path would go through zero or more faces overlapping with $\Omega^+$, all of which must have been assigned the same label by \cref{algo:tree}. Conversely, if for two points $x, y \in \Omega^+$ the evaluation of \cref{algo:eval} yields $\chi(x) = \chi(y)$, then there must be a sequence of leaf cells $U_1, U_2, \ldots, U_m$ with $x \in U_1$ and $y \in U_m$ such that, for each pair, $(U_i \cap U_{i+1}) \cap \{\phi > 0\}$ is nonempty. Therefore, there is a sequence of points $z_i \in (U_i \cap U_{i+1}) \cap \{\phi > 0\}$, and, because each $U_i \cap \Omega^+$ has exactly one component, there exists a path connecting $x \rightarrow z_1 \rightarrow z_2 \rightarrow \cdots \rightarrow z_{m-1} \rightarrow y$, and so $x$ and $y$ are path connected. The same arguments apply analogously to the negative region, and it is trivial to see that $\chi(x) \neq \chi(y)$ whenever $\phi(x)\phi(y) < 0$, since the negative and positive regions do not have overlapping label identifiers. In summary, when all the leaf cells satisfy the simply connected property, the output of \cref{algo:eval} is exact in the sense that $\chi(x) = \chi(y)$ if and only if $x$ and $y$ are path connected. In fact, we will observe in the next section that in many instances, the labeling algorithm continues to be exact even when the maximum recursion depth is reached.

\section{Numerical Experiments}
\label{sec:results}

In this section, we assess the effectiveness of the connected component labeling algorithm on randomly generated geometry in 2D and 3D as well as on particular examples exhibiting cusps, self-intersections, and other kinds of singularities.

\subsection{Randomly generated geometry}
\label{sec:random}

We consider here a class of randomly generated polynomials whose corresponding implicitly defined geometry exhibits a variety of characteristics such as high-curvature pieces and almost-touching components. The random polynomials are defined through the orthonormal Legendre polynomials, as follows. Let $p_i$, $i = 0,1,2,3,$ denote the first four univariate Legendre polynomials relative to the interval $[-1,1]$, i.e.,
\[ p_0(x) = 1, \quad p_1(x) = x, \quad p_2(x) = \tfrac12(3x^2 - 1), \quad p_3(x) = \tfrac12(5x^3 - 3x). \]
In $d$ dimensions, define a tensor-product degree $(3,\ldots,3)$ polynomial $\phi: [-1,1]^d \to \R$ as follows:
\begin{equation} \label{eq:random} \phi(x) = \sum_{i \in \{0,1,2,3\}^d} c_i\, \omega^i \prod_{\ell = 1}^d p_{i_\ell}(x_\ell), \end{equation}
where $c_i \in \R$ is a given set of $4^d$ coefficients, one for each multi-index $i \in \{0,1,2,3\}^d$, and $\omega > 0$ is a given (fixed) parameter such that $\omega^i := \smash{\prod_{\ell=1}^d \omega^{i_\ell}}$. With this set up, the class of randomly generated geometry is defined by polynomials of the form \eqref{eq:random} such that the coefficients $c_i$ are randomly and independently drawn from the uniform distribution on $[-1,1]$. Depending on the value of $\omega$, the factor $\omega^i$ controls the smoothness of the polynomial $\phi$ by damping (or magnifying) the higher-order oscillatory modes of the Legendre basis. In particular, we have set $\omega = 0.5$, empirically chosen to give a reasonable spread between relatively mild to more complex geometry involving multiple components of high degree curvature; examples are given in the next set of figures.

As discussed in the prior section, the connected component labeling algorithm is exact whenever the subdivision process finishes without reaching the maximum recursion depth $\Lambda$. Consequently, we expect that as $\Lambda$ increases, the percentage of cases in which the algorithm is certifiably correct will also increase. To quantify this, we generate a fixed number of random polynomials of the form \eqref{eq:random}, convert them to the Bernstein basis, and invoke \cref{algo:tree} on each instance on the reference domain $\mathcal U = [-1,1]^d$. Three possible outcomes are tallied by the following quantities:
\begin{itemize}
\item Let $\rho_{\shorteq}$ denote the percentage of cases in which line \ref{line:nonsc} is not executed, i.e., every leaf cell is simply connected so that the labeling is certifiably exact.
\item Let $\rho_\checkmark$ denote the percentage of cases in which line \ref{line:nonsc} is executed (i.e., the algorithm is uncertain about the topology), yet the labeling produced is exact.
\item Let $\rho_\times = 100\% - (\rho_{\shorteq}) - (\rho_\checkmark)$ denote the percentage of remaining cases, i.e., those in which two or more components have been artificially-glued.
\end{itemize}
Note that the last two scenarios require a notion of the ground truth so as to determine whether the labeling is exact or not; to that end, a reference solution for each of the randomly generated cases is computed by applying the base algorithm with unlimited recursion depth.\footnote{In every computed instance, the ground truth solution terminated in finite time.}

\begin{figure}%
\centering\small\sffamily%
\includegraphics[scale=0.9]{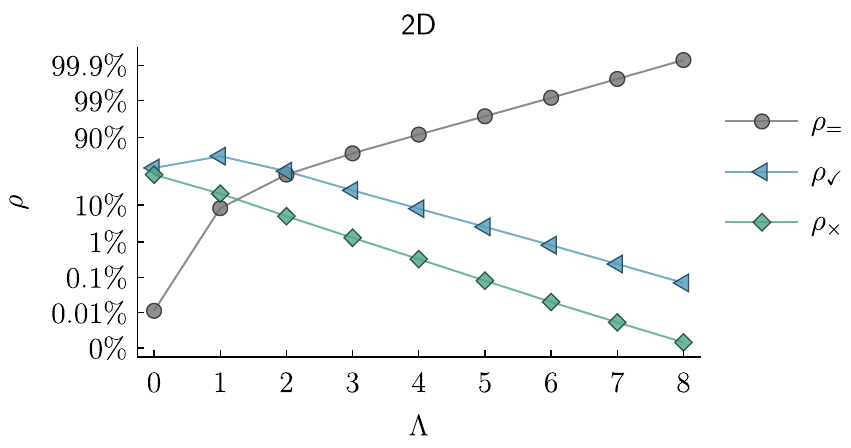}  \includegraphics[scale=0.9]{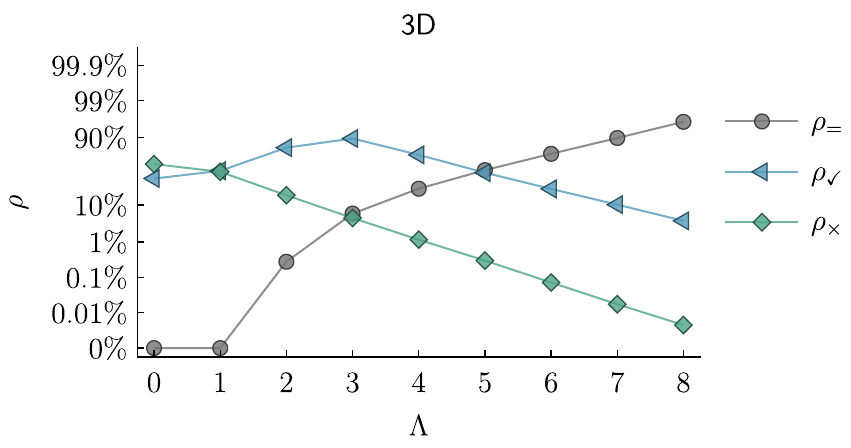}%
\caption{Percentage breakdown analysis corresponding to the three possible outcomes of the connected component labeling algorithm as it applied to the randomly generated geometry considered in \cref{sec:random}. As a function of the maximum recursion depth $\Lambda$, the quantities indicate the percentage of cases in which the algorithm: (i) yields certifiably correct labeling ($\rho_=$); (ii) is uncertain about the topology but produces correct labeling ($\rho_\checkmark$); and (iii) is uncertain about the topology such that two or more nearly-touching components are glued ($\rho_\times$); in particular, take note of the quasi-logarithmic scale.}
\label{fig:fractions}
\end{figure}

\cref{fig:fractions} plots $\rho_{\shorteq}$, $\rho_\checkmark$, and $\rho_\times$ as a function $\Lambda$.\footnote{The data shown in \cref{fig:fractions} was collected from ten million randomly generated cases in 2D and one million random examples in 3D, more than enough to stabilize the results.} On this class of randomly generated geometry, it is rare for the topology to be fully resolved with just one or two levels of subdivision (e.g., $\rho_{\shorteq}$ is less than 10\% when $\Lambda \leq 1$); however, after a few more levels of subdivision, $\rho_{\shorteq}$ increases to around 99\% when $\Lambda = 6$ in 2D, and around 90\% when $\Lambda = 7$ in 3D. In general, we observe that as $\Lambda$ increases, the number of cases in which the labeling algorithm is certifiably correct increases exponentially. Furthermore, note that $\rho_\times$ decays to zero at about the same exponential speed as $\rho_\checkmark$, but is a significantly small fraction thereof; as such, even when the algorithm reaches the maximum recursion depth and is consequently uncertain about the topology, the labeling continues to be exact in more cases than not. We also note that the failure rate $\rho_\times$ is somewhat similar in 2D and 3D, but the rate of certitude ($\rho_{\shorteq}$) in 3D is smaller than it is in the 2D setting. This can be attributed to a number of aspects, chief among which the 3D geometry is substantially more complex. Another major reason concerns the Bernstein polynomial range evaluation, which is generally less accurate in 3D versus 2D; consequently, a few more levels of subdivision are needed in order to pass the simply connected polynomial tests, and intuitively one may expect higher dimensions to require additional refinement power.

\newcommand{\lambdabox}[1]{\makebox[1.02in][c]{$\Lambda = #1$}}
\newcommand{\deltabox}[1]{\makebox[1.02in][c]{$\Delta = #1$}}

\begin{figure}%
\centering\small\sffamily%
\hspace{0.25mm}\includegraphics[width=6in]{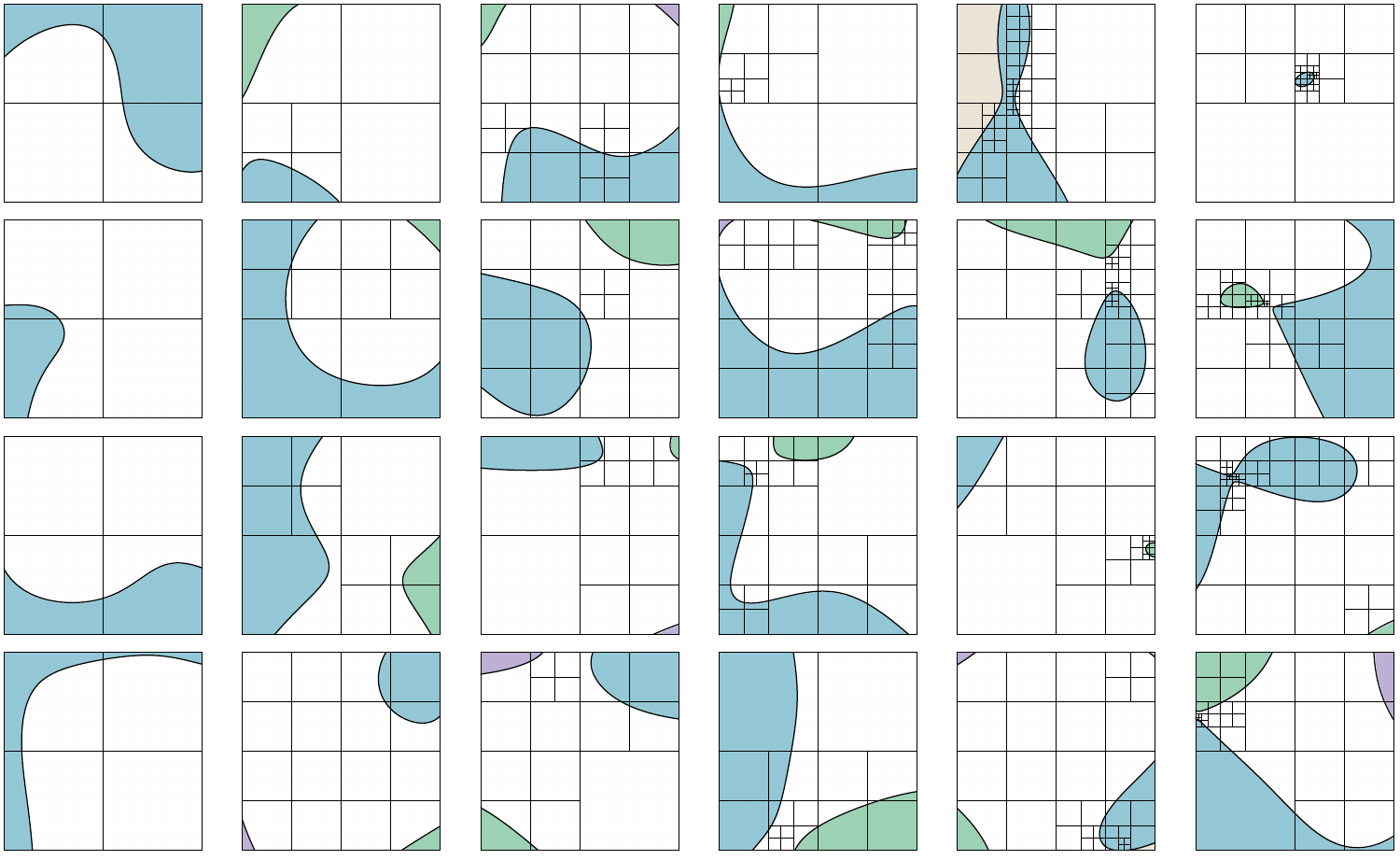}\\%
\deltabox{1}\deltabox{2}\deltabox{3}\deltabox{4}\deltabox{5}\deltabox{6}%
\caption{Examples of randomly generated geometry for which the connected component labeling algorithm yields certifiably correct results. The depth of the subdivision quadtree, $\Delta$, is indicated in each column, e.g., the far-left examples require just one subdivision step to be certain of the topology. Each example applies its own coloring scheme to illustrate the labeling output of the algorithm.}
\label{fig:a}
\end{figure}

\begin{figure}%
\centering\small\sffamily%
\hspace{0.25mm}\includegraphics[width=6in]{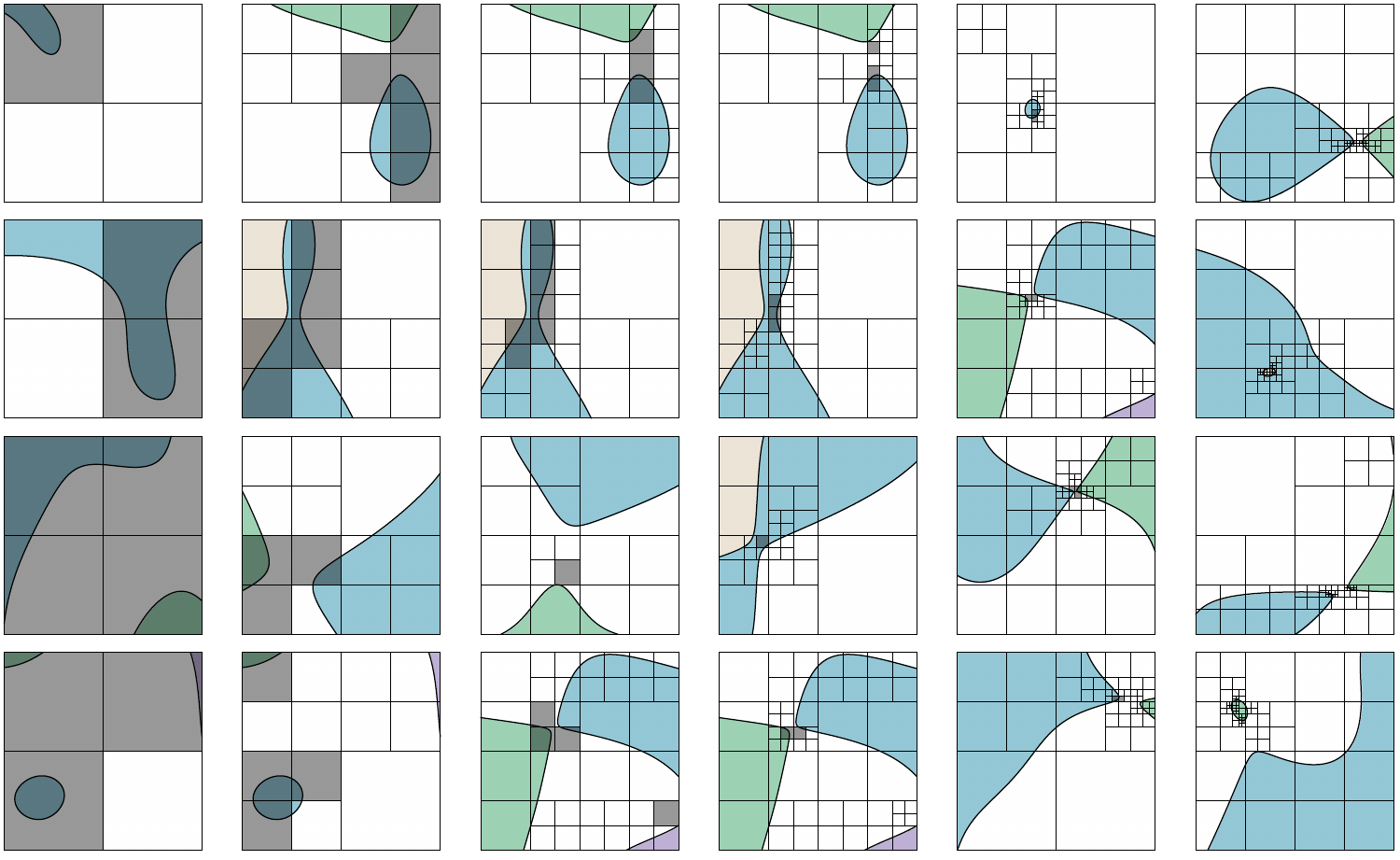}\\%
\lambdabox{1}\lambdabox{2}\lambdabox{3}\lambdabox{4}\lambdabox{5}\lambdabox{6}%
\caption{Examples of randomly generated geometry for which the connected component labeling algorithm reaches the maximum-imposed recursion depth $\Lambda$, is uncertain about the topology, yet correctly labels distinct components. Each example applies its own coloring scheme to illustrate the labeling output of the algorithm, while the shaded regions indicate the quadtree cells upon which the input polynomial is not simply connected, according to \cref{def:sc}.}
\label{fig:b}
\end{figure}

\begin{figure}%
\centering\small\sffamily%
\hspace{0.25mm}\includegraphics[width=6in]{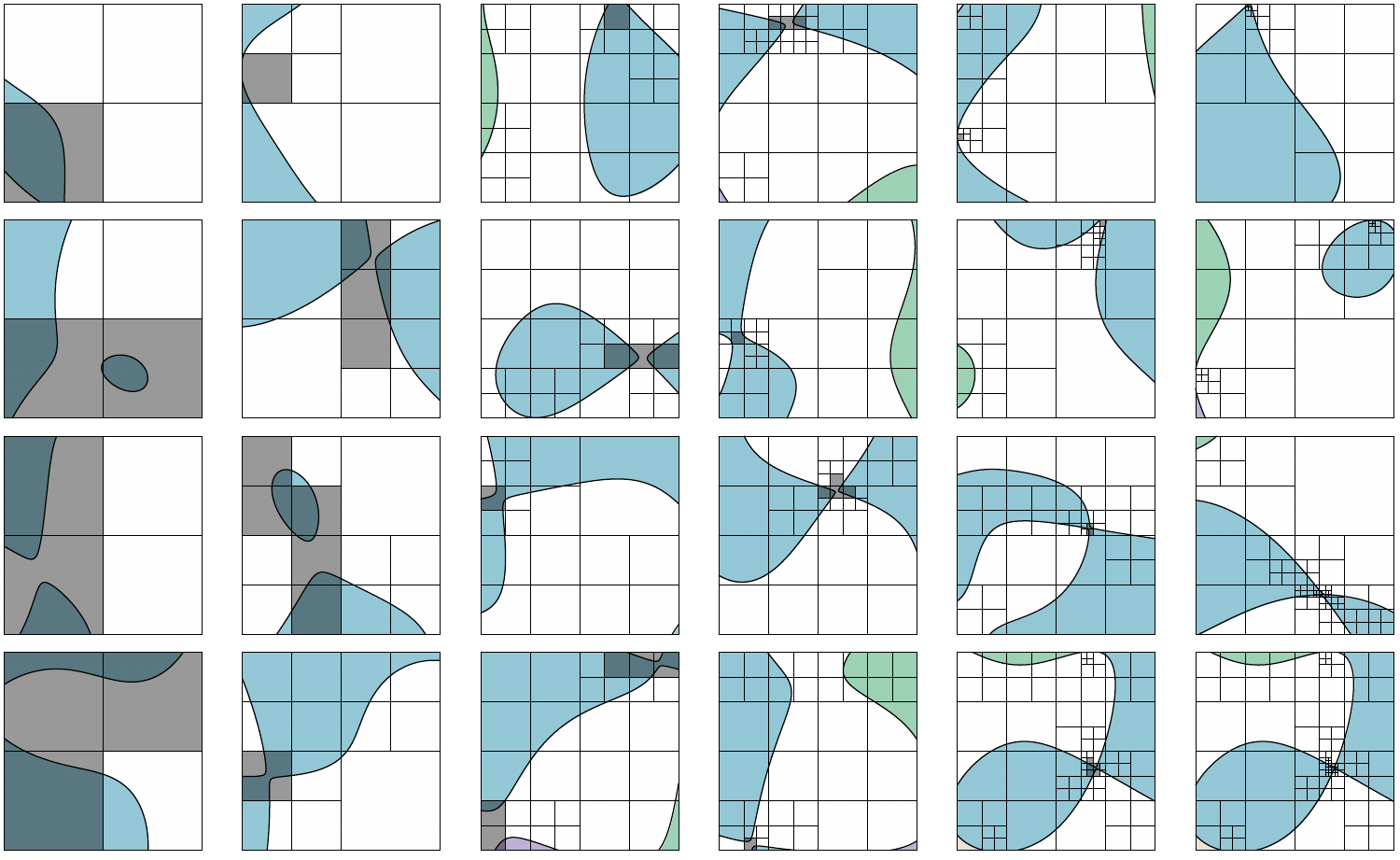}\\%
\lambdabox{1}\lambdabox{2}\lambdabox{3}\lambdabox{4}\lambdabox{5}\lambdabox{6}%
\caption{Examples of randomly generated geometry for which the connected component labeling algorithm reaches the maximum-imposed recursion depth $\Lambda$, is uncertain about the topology, and glues distinct components. Each example applies its own coloring scheme to illustrate the labeling of the algorithm; in particular, glued components have the same color. The shaded regions indicate the quadtree cells upon which the input polynomial is not simply connected, according to \cref{def:sc}; some of these cells are responsible for bridging the gap between glued components.}
\label{fig:c}
\end{figure}

\Cref{fig:a,fig:b,fig:c} provide some 2D examples demonstrating the above three scenarios. Each example illustrates the quadtree and the colors have a one-to-one correspondence with computed labeling; shaded cells correspond to those whereon $\phi$ is not simply connected, per the requirements of \cref{def:sc}. \cref{fig:a} shows some examples in which the algorithm is certain about the topology, i.e., every cell is simply connected; observe that, as expected, progressively more complex geometry can be resolved as further subdivision steps are permitted. \cref{fig:b} demonstrates cases in which the maximum recursion depth is reached (hence some cells are shaded), yet the labeling is still exact; these instances arise whenever a shaded cell $U$ contains at most one component each of $U \cap \{\phi > 0\}$ and $U \cap \{\phi < 0\}$, but for some reason the simply connected test did not hold. Often, this occurs when a saddle point or extrema is nearby the zero level set of $\phi$. Meanwhile, \cref{fig:c} demonstrates cases in which the subdivision is, in a sense, forced to stop early, leaving some components glued together. This scenario occurs whenever: (i) a shaded cell contains multiple components; and/or (ii) a shaded cell is $\cplus$-connected or $\cminus$-connected to a neighbor even when no path exists between them. Situation (ii) arises whenever the sign evaluation algorithm applied to the shared face is inexact in the sense that it predicts all possible signs, yet in actuality $\phi$ is uniformly signed on the face.\footnote{In 2D, on the randomly generated degree 3 polynomials under consideration, scenario (ii) never occurs because our particular range evaluation algorithm \cref{algo:range} on 1D edges is exact.}

\begin{figure}%
\centering\small\sffamily%
\includegraphics[width=6in]{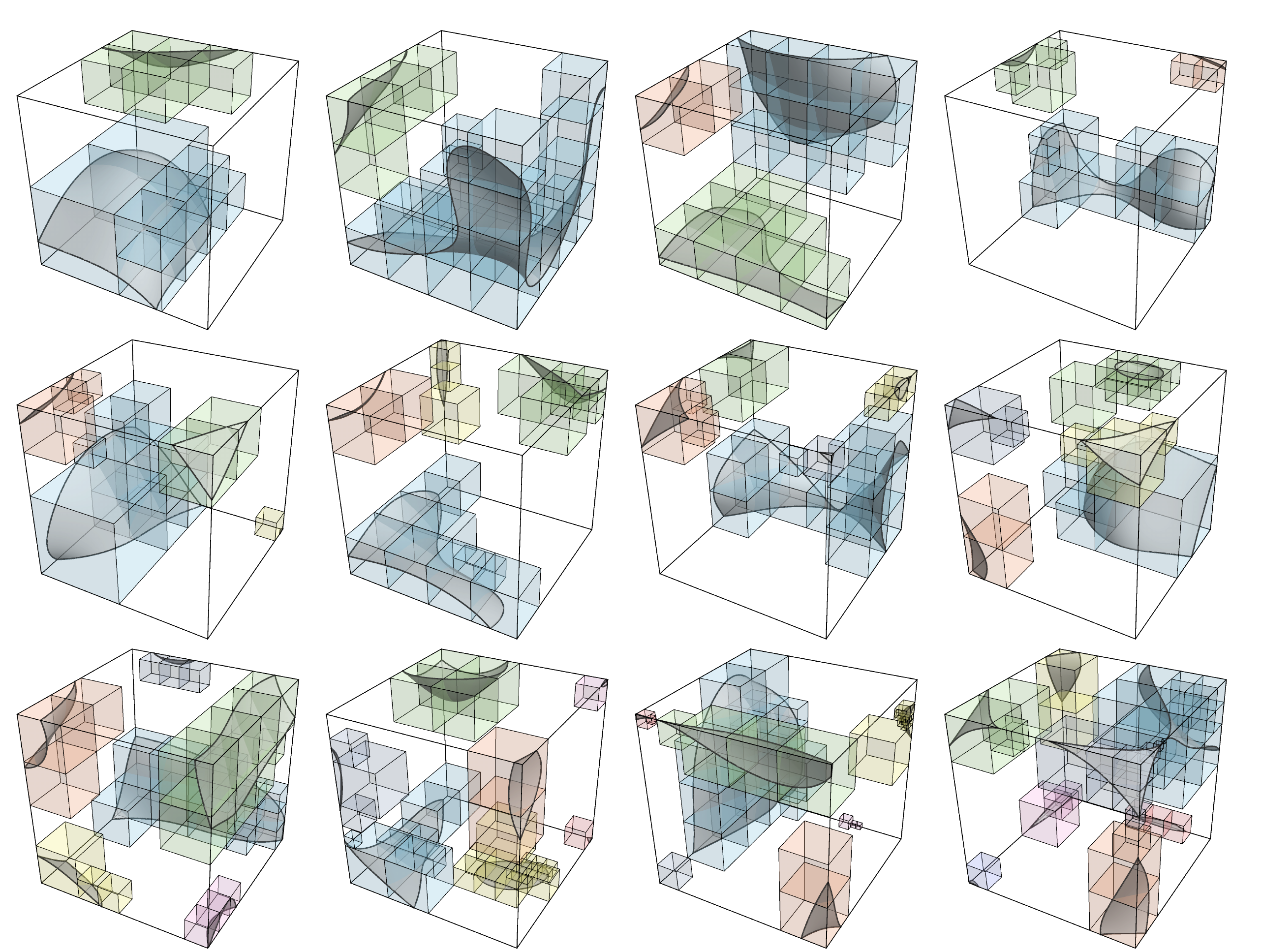}%
\caption{Examples illustrating the connected component labeling algorithm applied to randomly generated geometry in 3D. Only a subset of the octree is shown, corresponding to the cells for which $\Omega^+$ labels have been assigned. From top-left to bottom-right, the geometric complexity increases from two components up to eight.}
\label{fig:d}
\end{figure}

Similar behavior is obtained in 3D, though it is much less straightforward to visualize. \cref{fig:d} illustrates some examples. Note that in some cases, two cells with different labels can be touching, demonstrating that the sign evaluation algorithm was successful in disconnecting them. 

Intuitively, we expect the algorithm to glue components only when they are nearly touching. Indeed, the examples in \cref{fig:c} illustrate how the associated gaps are generally less than the size of the smallest cell, sometimes much smaller. To quantify this, we define an error metric which assigns a ``cost'' to gluing multiple components together. Roughly speaking, the cost is defined by the sum of the gaps between the agglomerated components. More precisely, given a set $V \subseteq {\mathcal U}$ and two points $x, y \in V$, we define the cost of connecting $x$ and $y$ by the shortest path connecting them, as measured by the arclength outside of $V$:
\[ {\mathcal C}(x, y; V) := \inf_{\gamma} \text{length} \bigl( \text{image}(\gamma) \setminus V \bigr), \]
where the infimum is taken over all paths $\gamma$ connecting $x$ and $y$ such that the part of $\gamma$ exterior to $V$ has measurable and finite length. The cost of gluing multiple connected components $U_1, \ldots, U_m$ is then defined as
\[ {\mathcal C}(U_1, \ldots, U_m) := \sup_{x, y \in V} {\mathcal C}(x, y; V) \text{ where } V = \bigcup_i U_i. \]
For example, the cost of gluing two connected components is the smallest straight-line gap between them. As another example, the cost of gluing the 1D sets $[0,1]$, $[2,3]$, and $[5,6]$ is $3$, being the sum of the two smaller gaps. 

\begin{figure}%
\centering\small\sffamily%
\includegraphics[scale=0.9]{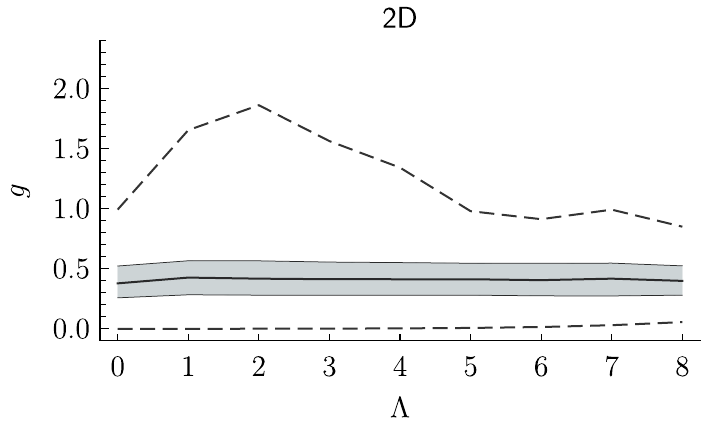}  \includegraphics[scale=0.9]{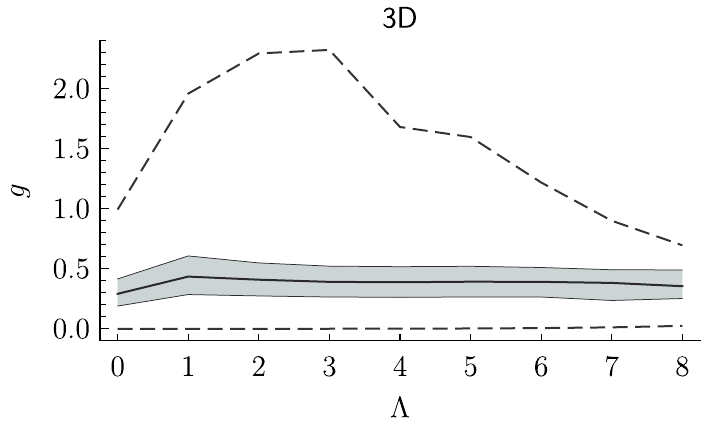}%
\caption{Gap distance between glued components, relative to the size of the smallest cell in the quadtree/octree. Aggregated over all relevant instances of the randomly generated geometry considered in \cref{sec:random}, i.e., those making up the $\rho_\times$ tally, the shaded region demarcates the first-to-third quartile spread of $g$, the interior solid line indicates the median, while the dashed lines indicate the minimum and maximum $g$ values.}
\label{fig:gaps}
\end{figure}

This metric is used to analyze the failure modes of the connected component labeling algorithm, as follows. For now, suppose $\Lambda$ is fixed. For each of the above-generated random polynomials $\phi_i$, we invoke \cref{algo:tree}. Using the reference solution, we determine which components have been glued together and compute the associated cost;\footnote{A separate algorithm has been developed to compute the gaps between the set of connected components of $\mathcal U \setminus \{\phi = 0\}$. The algorithm uses a recursive subdivision approach combined with a kind of adaptive point sampling. Its output is an interval bounding the exact gap value and can be made as tight as necessary; in particular, the gap calculations being presented here are sufficiently accurate for the analysis shown in \cref{fig:gaps}.} if no components were glued, the example is excluded from further analysis; if there are multiple sets of glued components, the cost of each set is measured and the maximal such cost is used. The output is the ``error'' $e_i$ for each example $i$ making up the $\rho_\times$ tally defined earlier; note that $e_i$ measures gaps of the input geometry, but its value depends on which components have been glued, which in turn depends on $\Lambda$. We expect $e_i$ to scale with the smallest cell size; in the present setting, for the reference domain $\mathcal U = [-1,1]^d$, the smallest cells have a diagonal length of $2^{1 - \Lambda} \sqrt{d}$. We therefore examine their ratio and define the overall error for polynomial $\phi_i$ by $g_i = e_i 2^{\Lambda - 1} / \sqrt{d}$. Aggregated over all such examples, we then compute the minimum and maximum value of $g_i$ as well as the median, first, and third quartiles. \cref{fig:gaps} displays the results corresponding to the set of 2D and 3D randomly generated polynomials. Once $\Lambda$ is sufficiently large so as to resolve the typical geometry (and for Bernstein subdivision to become sufficiently accurate), we observe that the gap between glued components is typically 25\%--50\% of the smallest cell. Meanwhile, the upper dashed curves of $g$ in \cref{fig:gaps} correspond to the rare cases in which two or more components are glued across a small chain of neighboring leaf cells. Notwithstanding this critical analysis, it is worth being reminded that, as \cref{fig:fractions} shows, the number of cases in which components are artificially glued decreases rapidly as $\Lambda$ is made larger.

\subsection{Additional examples}

\begin{figure}%
\centering\small\sffamily%
\begin{tabular}{ccc}
\includegraphics[width=1.7in]{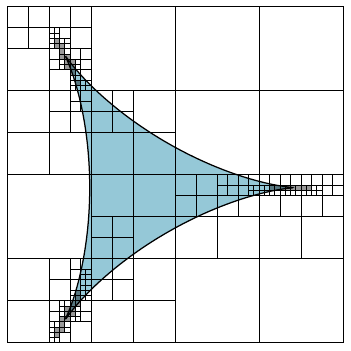} & \includegraphics[width=1.7in]{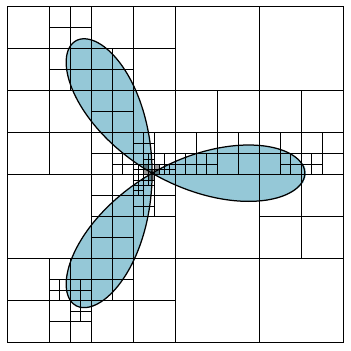} & \includegraphics[width=1.7in]{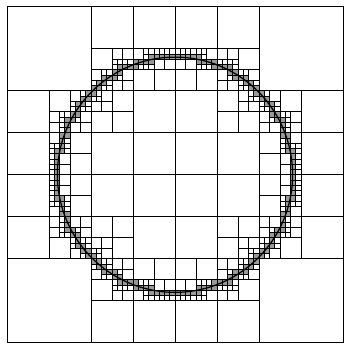} \\
(i) Deltoid, three cusps & (ii) Trifolium, self-intersection & (iii) Non-negative polynomial \\[1em]
\includegraphics[width=1.7in]{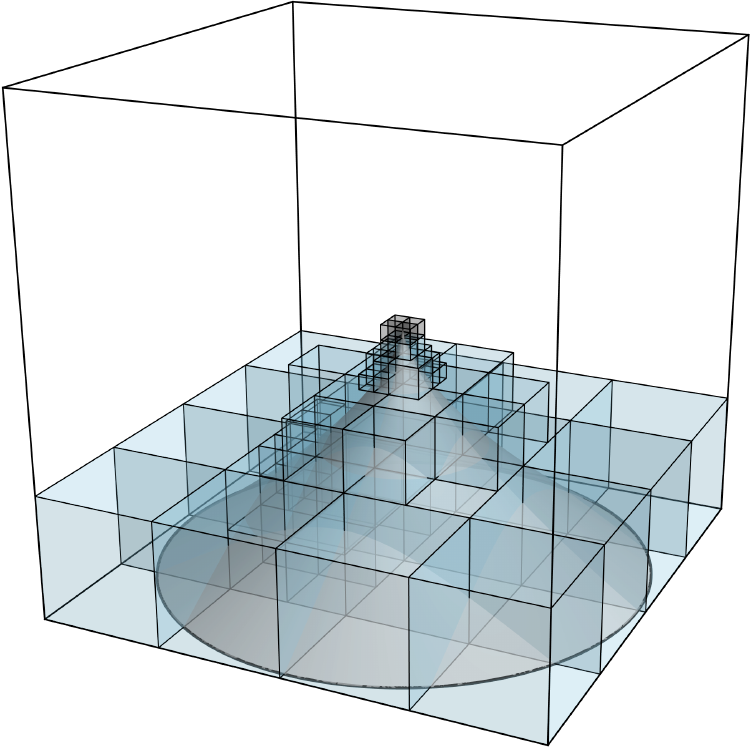} & \includegraphics[width=1.7in]{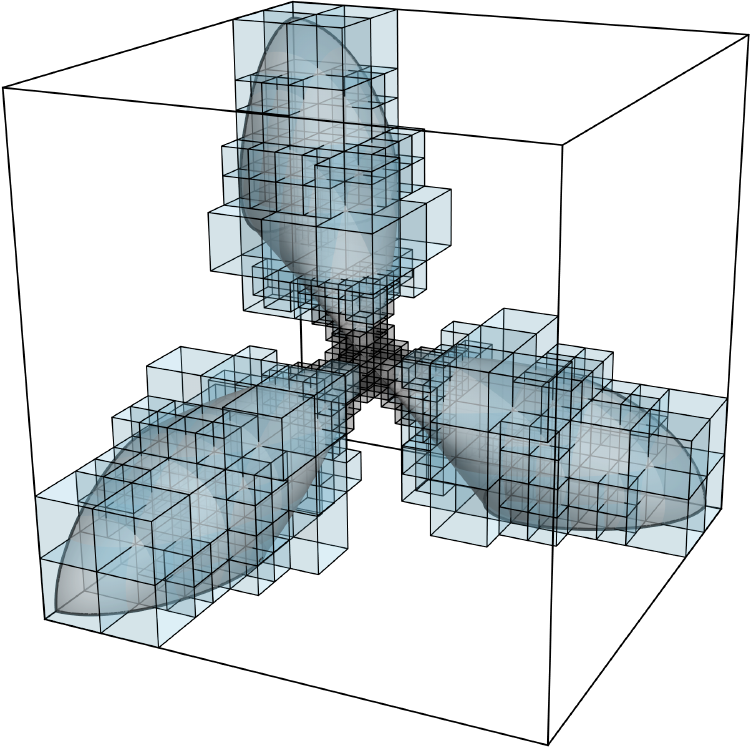} & \includegraphics[width=1.7in]{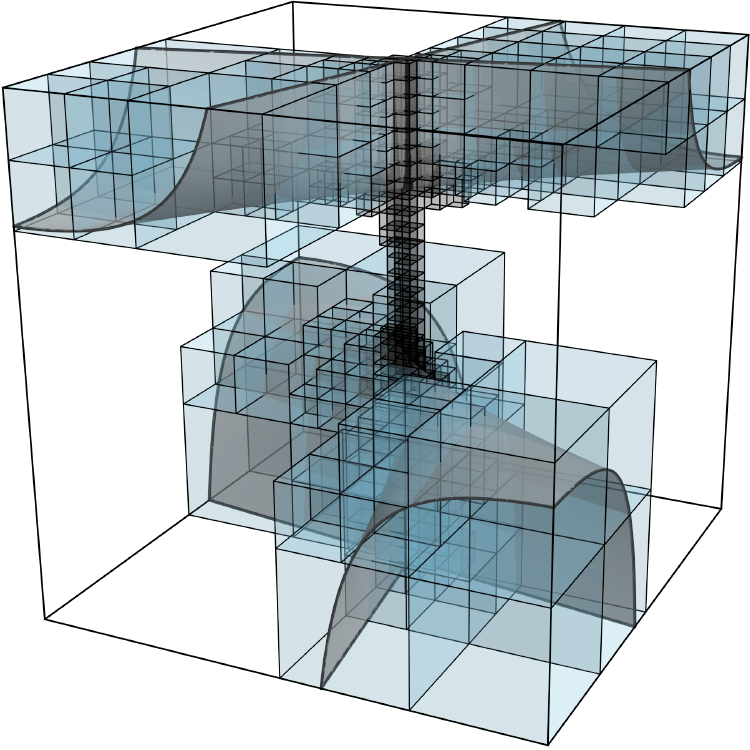} \\
(iv) Oloid, one cusp & (v) High-order junction & (vi) Line of singularities
\end{tabular}
\caption{Examples of test problems containing interfacial singularities.}
\label{fig:additional}
\end{figure}

\cref{fig:additional} illustrates some additional examples involving corners, cusps, self-intersections, and other kinds of singularities. These examples are more contrived in the sense they typically would not occur in a computational physics modeling problem, say; nonetheless, they serve to stress-test the connected component labeling algorithm in ways the randomly generated geometry did not. In particular, the examples demonstrate cases in which the implicitly-defined geometry contains points of singularity, i.e., points where both $\phi$ and its gradient are zero. In these cases, at least in comparison to the topology test of \cref{def:sc}, it is considerably more nuanced to determine whether the topology of $U \setminus \{\phi = 0\}$ is sufficiently simple on a subcell $U$. The examples of \cref{fig:additional} serve to illustrate the algorithm's behavior near singularities; in particular, (i) and (iv) show cases where the algorithm produces exact results despite the presence of singularities, while the remaining examples show cases where components are glued:
\begin{itemize}
\item \cref{fig:additional}(i) demonstrates a deltoid curve given by the zero level set of $\phi(x,y) = (x^2 + y^2)^2 + 18 (x^2 + y^2) - 8 (x^3 - 3 xy^2) - 27$; note how the subdivision focuses in on the singular portions of the curve situated at the deltoid's three cusps. In this case, the maximum-defined subdivision recursion depth is always met, but the resulting labeling is always exact.
\item \cref{fig:additional}(ii) demonstrates a trifolium curve, given by the zero level set $\phi(x,y) = (x^2 + y^2)^2 - x^3 + 3 x y^2$, containing a point of self-intersection at the origin. The subdivision process focuses in on the origin until it is forced to stop; the labeling algorithm subsequently glues the three connected components of $\{\phi > 0\}$, while the negative component $\{\phi < 0\}$ is unaffected.
\item Besides corners, cusps, and self-intersections, another kind of edge case concerns non-negative polynomials $\phi$ such that $\{\phi > 0\}$ contains multiple components. \cref{fig:additional}(iii) illustrates a case involving circular geometry, where $\phi(x,y) = (x^2 + y^2 - r^2)^2$. A ring of smallest-permitted subcells is obtained such that the components on either side of the interface are glued, as indicated by the solid white in the figure.
\item \cref{fig:additional}(iv) demonstrates a 3D case involving a single isolated cusp, given by conic-like oloid surface $\{\phi = 0\}$ where $\phi(x,y,z) = x^2 + y^2 + z^3$. Similar to the example in \cref{fig:additional}(i), the polynomial will never be simply connected on subcells sufficiently close to the cusp, and so the algorithm continues to subdivide until the maximum recursion depth is met. The labeling, however, is exact in this case.
\item \cref{fig:additional}(v) is a 3D analogue of (ii), and is given by the zero level set of $\phi(x,y,z) = (x + y + z - 1)(-x - y + z - 1)(x - y - z - 1)(-x + y - z - 1) - 2(x^2 + y^2 + z^2 - 3)^2$. The full surface on $\R^3$ contains six globular ``arms'' meeting at four junction points; we focus on just one of these singularities in the figure. In fixed-precision arithmetic, roundoff errors could lead to inconsistencies in the treatment of this geometry; the design choices in the present work result in all three components being glued.
\item Finally, \cref{fig:additional}(vi) demonstrates a non-trivial example in which the set of singularities is one-dimensional. Here, $\phi(x,y,z) = xyz - x^2 - y^2$, whose corresponding implicitly-defined geometry resembles two pinched pieces of paper, one above the other, rotated and connected via a line of singularities situated along the whole $z$-axis. In this example, the algorithm subdivides until the maximum-permitted resolution is met, with all subcells near or containing the $z$-axis marked as not simply connected. These leaf cells then form a bridge from one folded paper to the other, gluing all four of their respective components.
\end{itemize}
As an addendum, we point out that a tailor-made connected component labeling could be designed to specifically handle the kinds of geometry shown in \cref{fig:additional}. For example, in all of these cases the polynomial degree is small enough that it would be straightforward to apply cylindrical algebraic decomposition methods from real algebraic geometry. Our main purpose for these test problems is to demonstrate the ``failure modes'' of the connected component labeling algorithm in these exceptional situations.

\section{Concluding Remarks}
\label{sec:conclusion}

The connected component labeling algorithm developed here operates by recursively subdividing the constraint hyperrectangular domain $\mathcal U$ into progressively smaller subcells until the topology of $\Omega := \mathcal U \setminus \{\phi = 0\}$ thereon is sufficiently simple. This is achieved through a ``simply connected'' test exploiting some useful properties of the Bernstein polynomial basis, namely its effective methods for polynomial range evaluation. If the subdivision process succeeds, then the labeling is certifiably exact, i.e., $\chi(x) = \chi(y)$ if and only if $x, y \in \Omega$ are path-connected. At its core, however, the connected component labeling problem is ill-posed in the sense that tiny changes in the coefficients of the input polynomial $\phi$ could alter the number and topology of the connected components. We treated that here by allowing distinct components to be glued together, but only when the topology is uncertain: gluing occurs mainly when two components are nearly touching relative to the length scale of the smallest-permitted subcell; gluing also occurs for various edge cases such as interfacial self-intersections or junctions, or grazing-type singularities such as non-negative polynomials for which $\{\phi > 0\}$ has multiple components.

We have yet to discuss the speed of the algorithm. In general, the computational complexity is difficult to precisely quantify, being intricately dependent on the input polynomial $\phi$. For a fixed spatial dimension $d$ and bounded polynomial degree, the construction of the subdivision tree (\cref{algo:tree}) has a worst-case complexity $o(2^{d \Lambda})$; in practice, however, it is usually \textit{significantly} faster. Computation of $\chi(x)$ for some $x \in \Omega$ (\cref{algo:eval}) has a worst-case complexity ${\mathcal O}(\Lambda)$; in practice, however, it is usually \textit{significantly} faster. Indeed, for a given fixed polynomial, as soon as the subdivision process resolves its implicitly-defined geometry, the algorithmic complexity becomes independent of $\Lambda$.\footnote{In pathological cases, such as those illustrated in \cref{fig:additional}, the computational complexity depends on the measure of the input polynomial's set of singularities; in particular, if its set of singularities (inside $\mathcal U$) is nonempty and has maximal manifold dimension $d_s < d$, the worst-case complexity of \cref{algo:tree} is ${\mathcal O}(\lambda)$ when $d_s = 0$, or ${\mathcal O}(2^{d_s \Lambda})$ when $d_s > 0$, as $\Lambda \to \infty$.} A comparison to the high-order quadrature algorithms developed in \cite{SayeSP} is perhaps useful: solving the connected component labeling problem is about as fast as building a moderate order quadrature scheme. These quadrature algorithms are already ``fast'' in some sense, so the speed of the labeling algorithm is sufficient for our intended usage. As an additional indication, the cost of tree construction for the randomly generated geometry problems in \cref{sec:random} is about $8$ microseconds per instance in 2D and $0.2$ milliseconds per instance in 3D; about two thirds of that time is spent in the first phase of \cref{algo:tree}, about one third in the second phase, and negligible time in third phase. After tree construction and on the same set of test problems, the cost of evaluating $\chi(x)$ for a random point $x$ is about $0.05$ (resp., $0.1$) microseconds per evaluation in 2D (resp., 3D).\footnote{Timing measurements obtained on an Intel Xeon E3-1535m v6 laptop, single core, operating at approximately 3.5\,Ghz.}

Finally, we mention here some possibilities for extending the algorithmic approach. One immediate possibility is to develop a more sophisticated subcell topology test, e.g., $\Omega^\pm \cap U$ being star-connected or via an analysis of convexity \cite{LASSERRE2010912}. Accompanying this change, a more sophisticated sign evaluation algorithm might also be required, so as to appropriately path-connect neighboring cells, consistent with the updated topology test. Depending on the end application, however, these alterations might come with an increased computational cost. Another possibility is to replace the quadtree/octree-style subdivision with a k-d tree subdivision, or, more generally, a binary space partitioning (BSP) tree. In this setting, the most efficient subdivision process would likely occur by orienting to the features of the polynomial, e.g., by choosing hyperplanes which are locally parallel to its level sets. One would then have to generalize the range evaluation algorithms to polytopes and contend with more complex algorithms for building the cell connectivity. On the other hand, the quadtree/octree approach allows for a considerably simpler implementation, partly due to the use of tensor-product polynomials as well as dimension reduction and recursion techniques. As a final possibility, we note that the basic idea of the subdivision algorithm could be extended to non-polynomial level set functions; for example, one approach might be to apply automatic differentiation and interval arithmetic techniques (see, e.g., \cite{HighorderImplicitQuad}) to devise a suitable subcell topology test for arbitrary level set functions.

\section*{Acknowledgements}

This research was supported in part by the Applied Mathematics Program of the U.S.~Department of Energy (DOE) Office of Advanced Scientific Computing Research under contract number DE-AC02-05CH11231, and by a DOE Office of Science Early Career Research Program award.

\bibliographystyle{siamplain}
\bibliography{references}

\end{document}